\documentclass[12pt]{article}
\usepackage{amsmath, amssymb}

\newtheorem{theorem}{\bf Theorem}
\newtheorem{proposition}[theorem]{\bf Proposition}
\newtheorem{lemma}[theorem]{\bf Lemma}

\newtheorem{hypotheses}[theorem]{\bf Hypotheses}

\newcommand{\sect}[1]{\setcounter{equation}{0}\section{#1}}

\def\epsilon{\varepsilon}

\begin{document}

\Large \noindent 
{\bf Existence of traveling wave solutions
for a nonlocal monostable equation:
an abstract approach}

\vspace*{0.8em}

\normalsize
\noindent Hiroki Yagisita

\noindent
Department of Mathematics, 
Faculty of Science, 
Kyoto Sangyo University

\noindent
Motoyama, Kamigamo, Kita-Ku, Kyoto-City, 603-8555, Japan



\vspace*{1.6em}

\noindent {\bf Abstract} \ 
We consider the nonlocal analogue of the Fisher-KPP equation 
\[u_t=\mu*u-u+f(u),\]
where $\mu$ is a Borel-measure on $\mathbb R$ with $\mu(\mathbb R)=1$ 
and $f$ satisfies $f(0)=f(1)=0$ and $f>0$ in $(0,1)$. 
The equation may have a standing wave solution (a traveling wave solution with speed $0$) 
whose profile is a monotone but discontinuous function. 
We show that there is a constant $c_*$ 
such that it has a traveling wave solution with monotone profile and speed $c$ 
when $c\geq c_*$ 
while no periodic traveling wave solution with average speed $c$ when $c<c_*$. 
In order to prove it, we modify an abstract method for monotone semiflows by Weinberger. 
We note that the semiflow generated by the equation does not have compactness 
with respect to the compact-open topology. At the end of this paper, 
we propose a discrete Schrodinger model that describes the measurement process. 

\vspace*{0.4em}


\noindent 
equation, discrete monostable equation, 
nonlocal evolution equation, 

\noindent 
Fisher-Kolmogorov equation, multi-species mixture. 

\vspace*{0.8em}


\vspace*{0.8em}


\sect{Introduction}
\( \, \, \, \, \, \, \, \) 
In 1930, Fisher [6] introduced the reaction-diffusion equation $u_t=u_{xx}+u(1-u)$ 
as a model for the spatial spread of an advantageous form of a single gene 
in a population. He [7] found that there is a constant $c_*$ 
such that the equation has a traveling wave solution with speed $c$ when $c\geq c_*$ 
while it does no such solution when $c<c_*$. Kolmogorov, Petrovsky and Piskunov [14] 
investigated asymptotic behavior in the model. Since the pioneering works, there have been 
extensive studies on traveling waves and asymptotic behavior for monostable evolution systems.
In this paper, we consider the following nonlocal analogue of the Fisher-KPP equation: 
\[u_t=\mu*u-u+f(u).\]
Here, $\mu$ is a Borel-measure on $\mathbb R$ with $\mu(\mathbb R)=1$ and the convolution 
is defined by 
\[(\mu*u)(x)=\int_{y\in\mathbb R}u(x-y)d\mu(y)\]
for a bounded and Borel-measurable function $u$ on $\mathbb R$. The nonlinearity $f$ 
is a Lipschitz continuous function with $f(0)=f(1)=0$ and $f>0$ in $(0,1)$. Then, 
we would show that there is a constant $c_*$ 
such that the nonlocal monostable equation has a traveling wave solution 
with monotone profile and speed $c$ when $c\geq c_*$ 
while it does no periodic traveling wave solution with average speed c when $c<c_*$, 
if there is a positive constant $\lambda$ 
satisfying $\int_{y\in \mathbb R}e^{\lambda|y|}d\mu(y)<+\infty$. 
Further, we would also show that there is a smooth and monostable nonlinearity $f$ 
such that the equation has a standing wave solution (a traveling wave solution with speed $0$) 
whose profile is a monotone but discontinuous function, 
if $\mu$ satisfies the extra condition $\int_{y\in \mathbb R}yd\mu(y)>0$. 

For the nonlocal monostable equation, Schumacher [18, 19] proved that there is the minimal 
speed $c_*$ and the equation has a traveling wave solution with speed $c$ when $c\geq c_*$, 
if the nonlinearity $f$ satisfies the extra condition 
\[f(u)\leq f^\prime(0)u.\]   
On the other hand, Coville and Dupaigne [5] proved that the minimal speed $c_*$ is positive, 
the equation has a traveling wave solution with speed $c$ when $c\geq c_*$ 
and the profile of the solution is a smooth function, 
if the Borel-measure $\mu$ satisfies the extra condition
\[\mu((-\infty,-y))\equiv\mu((+y,+\infty)).\]
Further, Schumacher [18, 19] and Carr and Chmaj [2] studied uniqueness 
of traveling wave solutions for the equation. 
See, e.g., [4, 8, 9, 10, 11, 12, 13, 16, 20, 21, 23, 24] 
on traveling waves in various monostable evolution systems, 
[1, 3] nonlocal bistable equations and [17] Euler equation. 

In Section 2, we give abstract conditions such that a semiflow satisfying the conditions 
has a traveling wave solution with speed $c$ when $c\geq c_*$
while it does no periodic traveling wave solution with average speed c when $c<c_*$. 
We also note that it may not be compact with respect to the compact-open topology. 
In Section 3, we use idea in Weinberger [22] and Li, Weinberger and Lewis [15] 
to prove abstract theorems mentioned in Section 2. 
In Section 4, we precisely state our main results for the nonlocal monostable equation. 
In Section 5, we show that the semiflow generated by the nonlocal monostable equation 
satisfies the conditions given in Section 2 to prove the main results.  
At the end of this paper, we propose a discrete Schrodinger model 
that describes the measurement process. 
Another interest: validity of Boltzmann equation for multi-species mixture
(and/or semi-classical  gas). 
   
\sect{Abstract theorems for monotone semiflows}
\( \, \, \, \, \, \, \, \) 
In the abstract, we would treat a monostable evolution system. Put a set 
of functions on $\mathbb R$;
\[\mathcal M:=\{u\, |\, u \text{ is a monotone nondecreasing}\]
\[\text{ and left continuous function on } \mathbb R \text{ with } 0\leq u\leq 1\}.\] 

The followings are our basic conditions for discrete dynamical systems: 
\begin{hypotheses} \ 
Let $Q_0$ be a map from $\mathcal M$ into $\mathcal M$. 

{\rm (i)} \ $Q_0$ is continuous in the following sense: 
If a sequence $\{u_k\}_{k\in \mathbb N}\subset \mathcal M$ 
converges to $u\in \mathcal M$ uniformly on every bounded interval, 
then the sequence $\{Q_0[u_k]\}_{k\in \mathbb N}$ converges to $Q_0[u]$ almost everywhere. 

{\rm (ii)} \ $Q_0$ is order preserving; i.e., 
\[u_1\leq u_2 \Longrightarrow Q_0[u_1]\leq Q_0[u_2]\] 
for all $u_1$ and $u_2\in \mathcal M$. 

{\rm (iii)} \ $Q_0$ is translation invariant; i.e., 
\[(Q_0[u(\cdot)])(\cdot -x_0)=(Q_0[u(\cdot -x_0)])(\cdot)\] 
for all $u\in \mathcal M$ and $x_0\in \mathbb R$. 

{\rm (iv)} \ $Q_0$ is monostable; i.e., 
\[0<\alpha<1 \Longrightarrow \alpha < Q_0[\alpha]\]
for all constants $\alpha$. 
\end{hypotheses}
We note that the semiflow generated by a map $Q_0$ satisfying Hypotheses 1 
{\it may not be compact with respect to the compact-open topology}. 

\vspace*{0.4em} 

The following states that existence of suitable {\it super}-solutions 
of the form $\{v_n(x+cn)\}_{n=0}^\infty$ implies that of traveling wave solutions 
with speed $c$ in the discrete dynamical systems on $\mathcal M$: 
\begin{proposition} \ 
Let a map $Q_0:\mathcal M \rightarrow \mathcal M$ satisfy {\rm Hypotheses 1}, 
and $c\in \mathbb R$.  Suppose there exists a sequence $\{v_n\}^\infty_{n=0}\subset \mathcal M$ 
with $(Q_0[v_n])(x-c)\leq v_{n+1}(x)$, $\inf_{n=0, 1, 2, \cdots}{v_n(x)}\not\equiv 0$ 
and $\liminf_{n\rightarrow \infty}{v_n(x)}\not\equiv 1$.  Then, there exists $\psi\in \mathcal M$ 
with $(Q_0[\psi])(x-c)\equiv \psi(x)$, $\psi(-\infty)=0$ and $\psi(+\infty)=1$. 
\end{proposition}

\vspace*{0.2em} 

In the discrete dynamical system on $\mathcal M$ generated by a map $Q_0$ satisfying Hypotheses 1, 
if there is a {\it periodic} traveling wave {\it super-}solution with {\it average} speed $c$, 
then there is a traveling wave solution with speed $c$: 
\begin{theorem} \ 
Let a map $Q_0:\mathcal M \rightarrow \mathcal M$ satisfy {\rm Hypotheses 1}, 
and $c\in \mathbb R$.  Suppose there exist $\tau\in\mathbb N$ and $\phi\in\mathcal M$ 
with $({Q_0}^\tau[\phi])(x-c\tau)\leq \phi(x)$, $\phi\not\equiv 0$ and $\phi\not\equiv 1$. 
Then, there exists $\psi\in\mathcal M$ with  $(Q_0[\psi])(x-c)\equiv \psi(x)$, 
$\psi(-\infty)=0$ and $\psi(+\infty)=1$. 
\end{theorem}

\vspace*{0.2em} 

The infimum $c_*$ of the speeds of traveling wave solutions is not $-\infty$, 
and there is a traveling wave solution with speed $c$ when $c\geq c_*$:  
\begin{theorem} \ 
Suppose a map $Q_0:\mathcal M \rightarrow \mathcal M$ satisfies {\rm Hypotheses 1}. 
Then, there exists $c_*\in (-\infty, +\infty]$ such that the following holds {\rm :} 

Let $c\in \mathbb R$. Then, there exists $\psi\in\mathcal M$ 
with $(Q_0[\psi])(x-c\tau)\equiv \psi(x)$, 
$\psi(-\infty)=0$ and $\psi(+\infty)=1$ if and only if $c\geq c_*$. 
\end{theorem}

\vspace*{0.4em}

We add the following conditions to Hypotheses 1 for continuous dynamical systems 
on $\mathcal M$: 
\begin{hypotheses} \ 
Let $Q^t$ be a map from $\mathcal M$ to $\mathcal M$ for $t\in[0,+\infty)$. 

{\rm (i)} \ $Q$ is a semigroup; i.e., $Q^t\circ Q^s=Q^{t+s}$ for all $t$ and $s\in [0,+\infty)$. 

{\rm (ii)} \ $Q$ is continuous in the following sense: 
Suppose a sequence $\{t_k\}_{k\in \mathbb N}\subset[0,+\infty)$ converges to $0$, 
and $u\in \mathcal M$. Then, the sequence $\{Q^{t_k}[u]\}_{k\in \mathbb N}$ 
converges to $u$ almost everywhere. 
\end{hypotheses}

\vspace*{0.2em} 

As we would have Theorems 3 and 4 for the discrete dynamical systems, 
we would do the following two for the continuous dynamical systems:  
\begin{theorem} \ 
Let $Q^t$ be a map from $\mathcal M$ to $\mathcal M$ for $t\in[0,+\infty)$. 
Suppose $Q^t$ satisfies {\rm Hypotheses 1} for all $t\in (0,+\infty)$, 
and $Q$ {\rm Hypotheses 5}. Then, the following holds {\rm :} 

Let $c\in\mathbb R$. Suppose there exist $\tau\in(0,+\infty)$ and $\phi\in\mathcal M$ 
with $(Q^\tau[\phi])(x-c\tau)\leq \phi(x)$, $\phi\not\equiv 0$ and $\phi\not\equiv 1$. 
Then, there exists $\psi\in\mathcal M$ with $\psi(-\infty)=0$ and $\psi(+\infty)=1$ 
such that $(Q^t[\psi])(x-ct)\equiv \psi(x)$ holds for all $t\in[0,+\infty)$. 
\end{theorem}
\begin{theorem} \ 
Let $Q^t$ be a map from $\mathcal M$ to $\mathcal M$ for $t\in[0,+\infty)$. 
Suppose $Q^t$ satisfies {\rm Hypotheses 1} for all $t\in (0,+\infty)$, 
and $Q$ {\rm Hypotheses 5}. Then, there exists $c_*\in (-\infty, +\infty]$ 
such that the following holds {\rm :} 

Let $c\in \mathbb R$. Then, there exists $\psi\in\mathcal M$ 
with $\psi(-\infty)=0$ and $\psi(+\infty)=1$ 
such that $(Q^t[\psi])(x-ct)\equiv \psi(x)$ holds for all $t\in[0,+\infty)$
if and only if $c\geq c_*$. 
\end{theorem}

\sect{Proof of the abstract theorems}
\( \, \, \, \, \, \, \, \) 
In this section, we would modify the argument in Weinberger [22] 
and Li, Weinberger and Lewis [15] to prove the theorems stated in Section 2. 

\vspace*{0.2em} 

\begin{lemma} \ 
Let a sequence $\{u_k\}_{k\in \mathbb N}$ of monotone nondecreasing functions on $\mathbb R$ 
converge to a continuous function $u$ on $\mathbb R$ almost everywhere. 
Then, $\{u_k\}_{k\in \mathbb N}$ converges to $u$ uniformly on every bounded interval. 
\end{lemma}
{\bf Proof.} \ 
Let $C\in(0,+\infty)$ and $\varepsilon\in(0,+\infty)$. 
Then, there exists $\delta\in(0,+\infty)$ such that,
for any $y_1$ and $y_2\in[-C-1,+C+1]$, $|y_2-y_1|<\delta$ 
implies $|u(y_2)-u(y_1)|<\varepsilon/4$. So, we take 
$N\in\mathbb N$ and a sequence $\{x_n\}_{n=1}^N$ 
such that $\lim_{k\rightarrow\infty}u_k(x_n)=u(x_n)$, $-C-1\leq x_1\leq -C$, 
$x_n<x_{n+1}<x_n+\delta$ and $+C\leq x_N\leq +C+1$ hold. 

Let $k\in\mathbb N$ be sufficiently large. Then, 
$\max\{|u_k(x_n)-u(x_n)|\}_{n=1}^N<\epsilon/4$ holds. 
Let $x\in[-C,+C]$. There exists $n$ such that $x_n\leq x\leq x_{n+1}$ holds. 
So, we get $|u_k(x)-u(x)|\leq |u_k(x_n)-u(x)|+|u_k(x_{n+1})-u(x)|
\leq |u_k(x_n)-u(x_x)|+|u(x_n)-u(x)|
+|u_k(x_{n+1})-u(x_{n+1})|+|u(x_{n+1})-u(x)|<\varepsilon$. 
\hfill 
$\blacksquare$

\vspace*{0.8em}

The set of discontinuous points of a monotone function on $\mathbb R$ 
is at most countable. So, if a sequence $\{u_k\}_{k\in \mathbb N}$ 
of monotone functions on $\mathbb R$ 
converges to a monotone function $u$ on $\mathbb R$ 
at every continuous point of $u$, then  
it does almost everywhere. The converse also holds: 
\begin{lemma} \ 
Let a sequence $\{u_k\}_{k\in \mathbb N}$ of monotone nondecreasing functions on $\mathbb R$ 
converge to a monotone nondecreasing function $u$ on $\mathbb R$ almost everywhere. 
Then, $\lim_{k\rightarrow \infty}u_k(x)=u(x)$ holds 
for all continuous points $x\in\mathbb R$ of $u$. 
\end{lemma}
{\bf Proof.} \ 
We take $\underline{x}_n\in (x-2^{-n},x]$ and $\overline{x}_n\in [x,x+2^{-n})$ 
satisfying $\lim_{k\rightarrow\infty}
$ 
$
u_k(\underline{x}_n)=u(\underline{x}_n)$ 
and $\lim_{k\rightarrow\infty}u_k(\overline{x}_n)=u(\overline{x}_n)$ 
for $n\in\mathbb N$. Then, 
$u(\underline{x}_n) \leq \liminf_{k\rightarrow\infty}u_k(x) 
\leq \limsup_{k\rightarrow\infty}u_k(x) \leq u(\overline{x}_n)$ holds. 
Hence, we have $\lim_{k\rightarrow \infty}
$ 
$
u_k(x)=u(x)$ 
as $x$ is a continuous point of $u$. 
\hfill 
$\blacksquare$

\vspace*{0.8em}

Hypotheses 1 imply more continuity than Hypothesis 1 (i): 
\begin{proposition} \ 
Let a map $Q_0:\mathcal M \rightarrow \mathcal M$ 
satisfy {\rm Hypotheses 1} {\rm (i)}, {\rm (ii)} and {\rm (iii)}. 
Suppose a sequence $\{u_k\}_{k\in \mathbb N} \subset \mathcal M$ 
converges to $u\in \mathcal M$ almost everywhere. 
Then, $\lim_{k\rightarrow \infty}(Q_0[u_k])(x)=(Q_0[u])(x)$ holds 
for all continuous points $x\in\mathbb R$ of $Q_0[u]$. 
\end{proposition}
{\bf Proof.} \ 
We take a cut of function $\rho\in C^\infty(\mathbb R)$ with 
\[|x|\geq 1/2 \, \Longrightarrow \, \rho(x)=0,\]
\[|x|<1/2 \, \Longrightarrow \, \rho(x)>0\] 
and
\[\int_{x\in\mathbb R}\rho(x)dx=1.\] 
We put smooth functions 
\[\rho_n(\cdot):=2^n\rho(2^n\cdot),\]
\[\underline{u}^n(\cdot):=(\rho_n*u)(\cdot-2^{-(n+1)})\]
and 
\[\overline{u}^n(\cdot):=(\rho_n*u)(\cdot+2^{-(n+1)})\]
for $n\in\mathbb N$. Then, we obtain 
\[u(\cdot-2^{-n}) \leq \underline{u}^n(\cdot) 
\leq u(\cdot) \leq \overline{u}^n(\cdot) \leq u(\cdot+2^{-n}).\] 
Also, a sequence $\{\min\{u_k,\underline{u}^n\}\}_{k\in\mathbb N}$ 
converges to $\underline{u}^n$ almost everywhere 
and $\{\max\{u_k,\overline{u}^n\}\}_{k\in\mathbb N}$ $\overline{u}^n$. 
Hence, by Lemma 8, the sequence $\{\min\{u_k,\underline{u}^n\}\}_{k\in\mathbb N}$ converges 
to $\underline{u}^n$ uniformly on every bounded interval 
and $\{\max\{u_k,\overline{u}^n\}\}_{k\in\mathbb N}$ $\overline{u}^n$. 
Then, by Hypothesis 1 (i), the sequence $\{Q_0[\min\{u_k,\underline{u}^n\}]\}_{k\in\mathbb N}$ 
converges to $Q_0[\underline{u}^n]$ almost everywhere 
and $\{Q_0[\max\{u_k,\overline{u}^n\}]\}_{k\in\mathbb N}$ $Q_0[\overline{u}^n]$. 
From Hypothesis 1 (ii), $Q_0[\min\{u_k,\underline{u}^n\}] \leq Q_0[u_k] 
\leq Q_0[\max\{u_k,\overline{u}^n\}]$ also holds. 
Therefore, $Q_0[\underline{u}^n] \leq \liminf_{k\rightarrow\infty}Q_0[u_k] 
\leq \limsup_{k\rightarrow\infty}Q_0[u_k] \leq Q_0[\overline{u}^n]$ 
holds almost everywhere. So, by Hypotheses 1 (ii) and (iii), 
$Q_0[u](\cdot-2^{-n}) \leq \liminf_{k\rightarrow\infty}Q_0[u_k](\cdot) 
\leq \limsup_{k\rightarrow\infty}Q_0[u_k](\cdot) \leq Q_0[u](\cdot+2^{-n})$ 
holds almost everywhere. 
Hence, $\lim_{k\rightarrow\infty}Q_0[u_k](\cdot)=Q_0[u](\cdot)$ holds 
almost everywhere, because $\lim_{n\rightarrow\infty}Q_0[u](\cdot-2^{-n})
=\lim_{n\rightarrow\infty}Q_0[u](\cdot+2^{-n})=Q_0[u](\cdot)$ 
holds almost everywhere. So, from Lemma 9, 
$\lim_{k\rightarrow \infty}(Q_0[u_k])(x)=(Q_0[u])(x)$ holds 
for all continuous points $x\in\mathbb R$ of $Q_0[u]$.
\hfill 
$\blacksquare$

\vspace*{0.8em} 

Combining Proposition 10 with Helly's theorem, we can make 
the argument in Weinberger [22] and Li, Weinberger and Lewis [15] work 
to prove Proposition 2.

\vspace*{0.4em} 

\noindent
{\bf Proof of Proposition 2.}

We put $w(\cdot):=\lim_{h\downarrow +0}\inf_{n=0, 1, 2, \cdots}v_n(\cdot-h)$, 
and $u_0^k:=2^{-k}w\in\mathcal M$ for $k\in\mathbb N$. 
We also take functions $u_n^k\in\mathcal M$ such that 
\begin{equation}u_n^k(\cdot)=\max\{Q_0[u_{n-1}^k](\cdot-c),2^{-k}w(\cdot)\}\end{equation} 
holds for $k$ and $n\in\mathbb N$. 

We show 
\begin{equation}u_n^k\leq u_{n+1}^k.\end{equation} 
We have $u_0^k\leq u_1^k$. As $u_{n-1}^k\leq u_n^k$ holds, we get 
$Q_0[u_{n-1}^k]\leq Q_0[u_n^k]$ and $u_n^k\leq u_{n+1}^k$. 
In virtue of (3.2), we put $u^k:=\lim_{n\rightarrow\infty}u_n^k\in\mathcal M$. 
Then, by (3.1) and Proposition 10, 
\begin{equation}u^k(\cdot)=\max\{Q_0[u^k](\cdot-c),2^{-k}w(\cdot)\}\end{equation} 
holds. Because $\lim_{m\rightarrow\infty}Q_0[u^k(\cdot+m)]=Q_0[u^k(+\infty)]$ 
holds from Proposition 10, we have 
\[u^k(+\infty)=\lim_{m\rightarrow\infty}\max\{Q_0[u^k](m-c),2^{-k}w(m)\}\]
\[=\lim_{m\rightarrow\infty}\max\{Q_0[u^k(\cdot+m)](-c),2^{-k}w(m)\}\]
\[=\max\{Q_0[u^k(+\infty)],2^{-k}w(+\infty)\}.\] 
Hence, $u^k(+\infty)\geq Q_0[u^k(+\infty)]$ 
and $u^k(+\infty)\geq 2^{-k}w(+\infty)>0$ hold. So, from Hypothesis 1 (iv), 
we obtain 
\begin{equation}u^k(+\infty)=1.\end{equation} 

We show 
\begin{equation}u_n^k\leq v_n.\end{equation} 
We get $u_0^k\leq w\leq v_0$. As $u_{n-1}^k\leq v_{n-1}$ holds, we have 
\[Q_0[u_{n-1}^k](\cdot -c)\leq Q_0[v_{n-1}](\cdot-c)\leq v_n(\cdot)\] 
and $u_n^k\leq v_n$ because of $2^{-k}w\leq w\leq v_n$. From (3.5), 
\begin{equation}u^k(-\infty)\leq 
\lim_{m\rightarrow\infty}\liminf_{n\rightarrow\infty}v_n(-m)<1\end{equation} 
holds. Also, $\lim_{m\rightarrow\infty}Q_0[u^k(\cdot-m)]=Q_0[u^k(-\infty)]$ 
holds from Proposition 10. Hence, by (3.3), we have 
\[u^k(-\infty)
=\lim_{m\rightarrow\infty}\max\{Q_0[u^k](-m-c),2^{-k}w(-m)\}
\geq Q_0[u^k(-\infty)].\] 
So, from Hypothesis 1 (iv) and (3.6), we obtain 
\begin{equation}u^k(-\infty)=0.\end{equation}

In virtue of (3.4) and (3.7), there exists $x_k$ such that 
$u^k(-x_k)\leq 1/2 \leq \lim_{h\downarrow+0}u^k(-x_k+h)$ 
for $k\in\mathbb N$. We put $\psi^k(\cdot):=u^k(\cdot-x_k)\in\mathcal M$. 
Then, we have 
\begin{equation}\psi^k(0)\leq 1/2 \leq \lim_{h\downarrow+0}\psi^k(h)\end{equation} 
and 
\begin{equation}\psi^k(\cdot)=\max\{Q_0[\psi^k](\cdot-c),2^{-k}w(\cdot-x_k)\}\end{equation} 
from (3.3). By Helly's theorem, there exist a subsequence 
$\{k(n)\}_{n\in\mathbb N}$ and $\psi\in\mathcal M$ such that 
$\lim_{n\rightarrow \infty}\psi^{k(n)}(x)=\psi(x)$ holds 
for all continuous points $x\in\mathbb R$ of $\psi$. So, 
from (3.8), (3.9) and Proposition 10, 
\begin{equation}\psi(0)\leq 1/2 \leq \lim_{h\downarrow+0}\psi(h)\end{equation} 
and 
\begin{equation}\psi(\cdot)=Q_0[\psi](\cdot-c)\end{equation} 
holds. Because $\psi(-\infty)=Q_0[\psi(-\infty)]$ and $\psi(+\infty)=Q_0[\psi(+\infty)]$ 
also hold by (3.11) and Proposition 10, from Hypothesis 1 (iv) and (3.10), we have 
$\psi(-\infty)=0$ and $\psi(+\infty)=1$. 
\hfill 
$\blacksquare$ 

\vspace*{0.8em} 

\noindent
{\bf Proof of Theorem 3.} 

We take functions $v_n\in\mathcal M$ for $n=0, 1, 2, \cdots$ such that 
\[v_{n+m\tau}=({Q_0}^n[\phi])(\cdot-cn)\] 
holds for all $n=0, 1, 2, \cdots, \tau-1$ and $m=0, 1, 2, \cdots$. 
Then, we see 
\begin{equation}v_{n+1}(\cdot)\geq Q_0[v_n](\cdot-c)\end{equation} 
and 
\begin{equation}\liminf_{n\rightarrow\infty}v_n
=\inf_{n=0, 1, 2, \cdots}v_n=\min_{n=0, 1, 2, \cdots, \tau-1}v_n.\end{equation} 
We show $v_n(+\infty)>0$. We have $v_0(+\infty)>0$. 
As $v_{n-1}(+\infty)>0$ holds, we get $v_n(+\infty)\geq Q_0[v_{n-1}(+\infty)]>0$ 
by (3.12), Proposition 10, Hypotheses 1 (ii) and (iv). 
Hence, because $\lim_{m\rightarrow\infty}\min_{n=0, 1, 2, \cdots, \tau-1}v_n(m)>0$ holds, 
from (3.13), we see $\inf_{n=0, 1, 2, \cdots}v_n\not\equiv 0$. 
Because $\min_{n=0, 1, 2, \cdots, \tau-1}v_n\leq\phi$ holds, 
by (3.13) and $\phi(-\infty)<1$, 
we have $\liminf_{n\rightarrow\infty}v_n\not\equiv 1$. Therefore, 
by Proposition 2, there exists $\psi\in\mathcal M$ 
with ${Q_0}[\psi](\cdot-c)=\psi(\cdot)$, 
$\psi(-\infty)=0$ and $\psi(+\infty)=1$. 

\hfill 
$\blacksquare$ 

\vspace*{0.4em} 

\begin{lemma} \ 
Let a sequence $\{u_k\}_{k\in \mathbb N}$ of monotone nondecreasing functions on $\mathbb R$ 
converge to a monotone nondecreasing function $u$ on $\mathbb R$ almost everywhere. 
Then, $\lim_{k\rightarrow \infty}u_k(x-x_k)=u(x)$ holds 
for all sequences $\{x_k\}_{k\in \mathbb N} \subset \mathbb R$ 
with $\lim_{k\rightarrow \infty}x_k=0$ and continuous points $x\in\mathbb R$ of $u$. 
\end{lemma}
{\bf Proof.} \ 
We put $y_n:=\sup_{k=n, n+1, n+2, \cdots}|x_k|$ for $n\in\mathbb N$. 
Then, $u_k(\cdot-y_n) \leq u_k(\cdot-x_k) \leq u_k(\cdot+y_n)$ holds 
when $k \geq n$. Hence, 
$u(\cdot-y_n) \leq \liminf_{k\rightarrow\infty}u_k(\cdot-x_k) 
\leq \limsup_{k\rightarrow\infty}u_k(\cdot-x_k) \leq u(\cdot+y_n)$ 
holds almost everywhere. So, $\lim_{k\rightarrow\infty}u_k(\cdot-x_k)=u(\cdot)$ 
holds almost everywhere, because $\lim_{n\rightarrow\infty}u(\cdot-y_n)
=\lim_{n\rightarrow\infty}u(\cdot+y_n)=u(\cdot)$ 
holds almost everywhere. Hence, from Lemma 9, 
$\lim_{k\rightarrow \infty}u_k(x-x_k)=u(x)$ holds 
for all continuous points $x\in\mathbb R$ of $u$. 
\hfill 
$\blacksquare$

\vspace*{0.8em}

\noindent
{\bf Proof of Theorem 4.} 

[Step 1] \ Let $c_*\in[-\infty,+\infty]$ be the infimum of $c\in\mathbb R$ 
such that there exists $\psi\in\mathcal M$ 
with $Q_0[\psi](\cdot-c)=\psi(\cdot)$, 
$\psi(-\infty)=0$ and $\psi(+\infty)=1$. Then, we have the following: 
{\it Let $c\in\mathbb R$. 
Then, there exists $\psi\in\mathcal M$ 
with $Q_0[\psi](\cdot-c)=\psi(\cdot)$, 
$\psi(-\infty)=0$ and $\psi(+\infty)=1$ only if $c\geq c_*$.} 

[Step 2] \ In this step, we show the following: 
{\it Let $c\in(c_*,+\infty)$. 
Then, there exists $\psi\in\mathcal M$ 
with $Q_0[\psi](\cdot-c)=\psi(\cdot)$, 
$\psi(-\infty)=0$ and $\psi(+\infty)=1$.} 

There exist $c_0\in(-\infty,c)$ and $\phi\in\mathcal M$ 
with $Q_0[\phi](\cdot-c_0)=\phi(\cdot)$, 
$\phi(-\infty)=0$ and $\phi(+\infty)=1$. 
Then, because we have $Q_0[\phi](\cdot-c)\leq\phi(\cdot)$, 
by Theorem 3, there exists $\psi\in\mathcal M$ 
with $Q_0[\psi](\cdot-c)=\psi(\cdot)$, 
$\psi(-\infty)=0$ and $\psi(+\infty)=1$. 

[Step 3] \ In this step, we show the following: 
{\it Let $c_*\in\mathbb R$. 
Then, there exists $\psi\in\mathcal M$ 
with $Q_0[\psi](\cdot-c_*)=\psi(\cdot)$, 
$\psi(-\infty)=0$ and $\psi(+\infty)=1$.} 

In virtue of Step 2, there exists $\psi_k\in\mathcal M$ 
with $Q_0[\psi_k](\cdot-(c_*+2^{-k}))=\psi_k(\cdot)$, 
$\psi_k(-\infty)=0$ and $\psi_k(+\infty)=1$ for $k\in\mathbb N$. 
We also take $x_k$ such that 
$\psi_k(-x_k)\leq 1/2 \leq \lim_{h\downarrow+0}\psi_k(-x_k+h)$, 
and put $\psi^k(\cdot):=\psi_k(\cdot-x_k)\in\mathcal M$. 
Then, we have 
\begin{equation}\psi^k(0)\leq 1/2 \leq \lim_{h\downarrow+0}\psi^k(h)\end{equation} 
and 
\begin{equation}Q_0[\psi^k(\cdot-2^{-k})](\cdot-c_*)=\psi^k(\cdot).\end{equation} 
By Helly's theorem, there exist a subsequence 
$\{k(n)\}_{n\in\mathbb N}$ and $\psi\in\mathcal M$ such that 
$\lim_{n\rightarrow \infty}\psi^{k(n)}(x)=\psi(x)$ holds 
for all continuous points $x\in\mathbb R$ of $\psi$. Also, by Lemma 11, 
$\lim_{n\rightarrow \infty}\psi^{k(n)}(x-2^{-k(n)})=\psi(x)$ holds 
for all continuous points $x\in\mathbb R$ of $\psi$. 
Therefore, from (3.14), (3.15) and Proposition 10, 
\begin{equation}\psi(0)\leq 1/2 \leq \lim_{h\downarrow+0}\psi(h)\end{equation} 
and 
\begin{equation}Q_0[\psi](\cdot-c_*)=\psi(\cdot)\end{equation} 
holds. Because $Q_0[\psi(-\infty)]=\psi(-\infty)$ 
and $Q_0[\psi(+\infty)]=\psi(+\infty)$ 
also hold by (3.17) and Proposition 10, from Hypothesis 1 (iv) and (3.16), 
we have $\psi(-\infty)=0$ and $\psi(+\infty)=1$.

[Step 4] \ Finally, we show $c_*\in(-\infty,+\infty]$. 

Suppose $c_*=-\infty$. Then, in virtue of Step 2, 
there exists $\phi_k\in\mathcal M$ 
with $Q_0[\phi_k](\cdot+2^k)=\phi_k(\cdot)$, 
$\phi_k(-\infty)=0$ and $\phi_k(+\infty)=1$ for $k\in\mathbb N$. 
We also take $x_k$ such that 
$\phi_k(-x_k)\leq 1/2 \leq \lim_{h\downarrow+0}\phi_k(-x_k+h)$, 
and put $\phi^k(\cdot):=\phi_k(\cdot-x_k)\in\mathcal M$. 
Then, we have 
\begin{equation}\phi^k(0)\leq 1/2 \leq \lim_{h\downarrow+0}\phi^k(h)\end{equation} 
and 
\begin{equation}Q_0[\phi^k(\cdot+2^{k})](\cdot)=\phi^k(\cdot).\end{equation} 
Put $\chi\in\mathcal M$ such that $\chi(x)=0 \, (x\leq 0)$ 
and $\chi(x)=1/2 \, (0<x)$. Then, $\chi\leq\phi^k$ holds from (3.18). 
Hence, by (3.18) and (3.19), we see $Q_0[\chi(\cdot+2^{k})](0)\leq 1/2$. 
So, from $\lim_{k\rightarrow\infty}\chi(\cdot+2^{k})=1/2$ and Proposition 10, 
we obtain $Q_0[1/2]\leq 1/2$. This is a contradiction with Hypothesis 1 (iv). 
\hfill 
$\blacksquare$

\vspace*{0.4em} 

\begin{lemma} \ 
Let $Q^t$ be a map from $\mathcal M$ to $\mathcal M$ for $t\in[0,+\infty)$. 
Suppose $Q$ satisfies {\rm Hypothesis 5 (ii)}. 
Then, $\lim_{t\rightarrow 0}(Q^t[u])(x-ct)=u(x)$ holds 
for all $c\in\mathbb R$, $u\in \mathcal M$ 
and continuous points $x\in\mathbb R$ of $u$. 
\end{lemma}
{\bf Proof.} \ 
Let a sequence $\{t_k\}_{k\in \mathbb N}\subset[0,+\infty)$ converge to $0$. 
Then, by Hypothesis 5 (ii) and Lemma 11, 
$\lim_{k\rightarrow \infty}Q^{t_k}[u](x-ct_k)=u(x)$ holds 
for all continuous points $x\in\mathbb R$ of $u$. 
\hfill 
$\blacksquare$

\vspace*{0.8em} 

\noindent 
{\bf Proof of Theorem 6.} 

By Theorem 3, there exists $\psi_k\in\mathcal M$ 
with $Q^{\frac{\tau}{2^k}}[\psi_k](\cdot-\frac{c\tau}{2^k})=\psi_k(\cdot)$, 
$\psi_k(-\infty)=0$ and $\psi_k(+\infty)=1$ for $k\in\mathbb N$. 
We also take $x_k$ such that 
$\psi_k(-x_k)\leq 1/2 \leq \lim_{h\downarrow+0}\psi_k(-x_k+h)$, 
and put $\psi^k(\cdot):=\psi_k(\cdot-x_k)\in\mathcal M$. 
Then, we have 
\begin{equation}\psi^k(0)\leq 1/2 \leq \lim_{h\downarrow+0}\psi^k(h)\end{equation} 
and 
\begin{equation}Q^{\frac{\tau}{2^k}}[\psi^k](\cdot-\frac{c\tau}{2^k})
=\psi^k(\cdot).\end{equation} 
By Helly's theorem, there exist a subsequence 
$\{k(n)\}_{n\in\mathbb N}$ and $\psi\in\mathcal M$ such that 
$\lim_{n\rightarrow \infty}\psi^{k(n)}(x)=\psi(x)$ holds 
for all continuous points $x\in\mathbb R$ of $\psi$. 

Let $k_0\in\mathbb N$ and $m_0\in\mathbb N$. 
As $n\in\mathbb N$ is sufficiently large, 
\[Q^{\frac{m_0\tau}{2^{k_0}}}[\psi^{k(n)}](\cdot-c\frac{m_0\tau}{2^{k_0}})\] 
\[=(Q^{\frac{\tau}{2^{k(n)}}})^{m_02^{k(n)-k_0}}
[\psi^{k(n)}](\cdot-\frac{c\tau}{2^{k(n)}}m_02^{k(n)-k_0})
=\psi^{k(n)}(\cdot)\] 
holds because of $k(n)\geq k_0$ and (3.21). Therefore, by Proposition 10, we obtain 
\begin{equation}Q^{\frac{m_0\tau}{2^{k_0}}}[\psi](\cdot-c\frac{m_0\tau}{2^{k_0}})
=\psi(\cdot).\end{equation} 
From (3.20), we also see 
\begin{equation}\psi(0)\leq 1/2 \leq \lim_{h\downarrow+0}\psi(h).\end{equation} 

Let $t\in[0,+\infty)$. Then, by (3.22), 
there exists a sequence $\{t_k\}_{k\in\mathbb N}\subset[0,+\infty)$ 
with $\lim_{k\rightarrow\infty}t_k=0$ such that 
$Q^{t+t_k}[\psi](\cdot-c(t+t_k))=\psi(\cdot)$ 
holds for all $k\in\mathbb N$. So, by 
$Q^{t_k}[Q^t[\psi](\cdot-ct)](\cdot-ct_k)=Q^{t+t_k}[\psi](\cdot-c(t+t_k))$ 
and Lemma 12, we obtain 
\[Q^t[\psi](\cdot-ct)=\psi(\cdot).\] 
Hence, because $Q^t[\psi(-\infty)]=\psi(-\infty)$ 
and $Q^t[\psi(+\infty)]=\psi(+\infty)$ hold by Proposition 10, 
from (3.23), we see $\psi(-\infty)=0$ and $\psi(+\infty)=1$. 
\hfill 
$\blacksquare$

\vspace*{0.8em} 

\noindent 
{\bf Proof of Theorem 7.} 

In virtue of Theorem 4, we take $c_*\in (-\infty, +\infty]$ 
such that the following holds: 
{\it Let $c\in \mathbb R$. Then, there exists $\phi\in\mathcal M$ 
with $(Q^1[\phi])(\cdot-c)\equiv \phi(\cdot)$, 
$\phi(-\infty)=0$ and $\phi(+\infty)=1$ if and only if $c\geq c_*$.}  

Then, from Theorem 6, we have the conclusion of this theorem. 
\hfill 
$\blacksquare$

\sect{The main results for the nonlocal monostable equation}
\( \, \, \, \, \, \, \, \) 
Let a Lipschitz continuous function $f$ on $\mathbb R$ be a monostable nonlinearity; 
$f(0)=f(1)=0$ and $f(u)>0$ in $(0,1)$. Let a Borel-measure $\mu$ on $\mathbb R$ 
satisfy $\mu(\mathbb R)=1$. Then, we consider the following nonlocal monostable equation: 
\begin{equation}u_t=\mu*u-u+f(u),\end{equation} 
where $(\mu*u)(x):=\int_{y\in\mathbb R}u(x-y)d\mu(y)$ 
for a bounded and Borel-measurable function $u$ on $\mathbb R$. Then, 
$G(u):=\mu*u-u+f(u)$ is a map from the Banach space $L^\infty(\mathbb R)$ 
into $L^\infty(\mathbb R)$ and it is Lipschitz continuous  
(we note that $u(x-y)$ is a Borel-measurable function on $\mathbb R^2$, 
and $\|u\|_{L^\infty(\mathbb R)}=0$ implies 
$\|\mu*u\|_{L^1(\mathbb R)}\leq 
\int_{y\in\mathbb R}(\int_{x\in\mathbb R}|u(x-y)|dx)d\mu(y)$=0). 
So, because the standard theory of ordinary differential equations works, 
we have well-posedness of (4.1) and the equation 
generates a flow in $L^\infty(\mathbb R)$. 
The following gives two positively invariant sets: 
\begin{proposition} \ 
If $u_0\in L^\infty(\mathbb R)$ satisfies $0\leq u_0\leq 1$, then 
there exists a solution $\{u(t)\}_{t\in[0,+\infty)}\subset L^\infty(\mathbb R)$ 
to (4.1) with $u(0)=u_0$ and $0\leq u(t)\leq 1$. 

For any $u_0\in \mathcal M$, then there exists 
a solution $\{u(t)\}_{t\in[0,+\infty)}\subset \mathcal M$ 
to (4.1) with $u(0)=u_0$. 
\end{proposition}

\vspace*{0.4em}

If the semiflow generated by (4.1) has a {\it periodic} traveling wave solution 
with {\it average} speed $c$ (even if the profile is not a monotone function), 
then it dose a traveling wave solution with {\it monotone} profile and speed $c$:  
\begin{theorem} \ 
Let a Borel-measure $\mu$ have $\lambda\in(0,+\infty)$ satisfying 
\begin{equation}\int_{y\in\mathbb R}e^{\lambda|y|}d\mu(y)
<+\infty,\end{equation} 
and $c\in\mathbb R$. Suppose there exist $\tau\in(0,+\infty)$ 
and a solution $\{u(t,x)\}_{t\in\mathbb R} 
\subset L^\infty(\mathbb R)$ to (4.1) with $0\leq u(t,x)\leq 1$, 
$\lim_{x\rightarrow+\infty} u(t,x)=1$ and $\|u(t,x)-1\|_{L^\infty(\mathbb R)}\not=0$ 
such that \[u(t+\tau,x)=u(t,x+c\tau)\] holds for all $t$ and $x\in\mathbb R$.  
Then, there exists $\psi\in\mathcal M$ with $\psi(-\infty)=0$ and $\psi(+\infty)=1$ 
such that $\{\psi(x+ct)\}_{t\in\mathbb R}$ is a solution to (4.1). 
\end{theorem}

\vspace*{0.4em} 

The infimum $c_*$ of the speeds of traveling wave solutions is not $\pm\infty$, 
and there is a traveling wave solution with speed $c$ when $c\geq c_*$: 
\begin{theorem} \ 
Let a Borel-measure $\mu$ have $\lambda\in(0,+\infty)$ satisfying (4.2). 
Then, there exists $c_*\in\mathbb R$ such that the following holds {\rm :} 

Let $c\in \mathbb R$. Then, there exists $\psi\in\mathcal M$ 
with $\psi(-\infty)=0$ and $\psi(+\infty)=1$ 
such that $\{\psi(x+ct)\}_{t\in\mathbb R}$ is a solution to (4.1) 
if and only if $c\geq c_*$. 
\end{theorem}

\vspace*{0.4em} 

The solutions to (4.1) are continuous in $L^\infty(\mathbb R)$. 
Hence, if the profile of a traveling wave solution with speed $c\not=0$ 
is monotone, then it is a continuous function on $\mathbb R$. 
However, for some nonlinearity $f$, 
if the profile of a standing wave solution (a traveling wave solution with speed $0$) 
is a monotone function, then it is a discontinuous one: 
\begin{proposition} \ 
Let a nonlinearity $f\in C^1(\mathbb R)$ satisfy $\max_{u\in[0,1]}f_u(u)>1$, 
and $\psi\in\mathcal M$ $\psi(-\infty)=0$ and $\psi(+\infty)=1$. 
Suppose $u(t,x):=\psi(x)$ is a solution to (4.1). 
Then, $\psi$ is a discontinuous function. 
\end{proposition}

\vspace*{0.4em} 

Coville and Dupaigne [5] shown that the minimal speed $c_*$ is positive, 
if the Borel-measure $\mu$ satisfies the following extra condition: 
\[\mu((-\infty,-y))\equiv\mu((+y,+\infty)).\] 
It implies $\int_{y\in\mathbb R}yd\mu(y)=0$. On the other hand, 
if the Borel-measure $\mu$ satisfies $\int_{y\in\mathbb R}yd\mu(y)>0$, 
then the minimal speed $c_*$ is negative and, so, 
the semiflow has a standing wave solution 
(a traveling wave solution with speed $0$) 
for a sufficiently small nonlinearity $f$: 
\begin{proposition} \ 
Let a Borel-measure $\mu$ have $\lambda\in(0,+\infty)$ satisfying (4.2), 
and $\int_{y\in\mathbb R}yd\mu(y)>0$. Then, there exists $\gamma\in(0,+\infty)$ 
such that the minimal speed $c_*$ is negative 
when $f(u)\leq \gamma u \ (0\leq u\leq 1)$. 
\end{proposition}

\sect{Semiflows generated by nonlocal monostable equations}
\( \, \, \, \, \, \, \, \) 
In this section, we show the results for the nonlocal monostable equation (4.1) 
stated in Section 4. 

First, we see a comparison theorem: 
\begin{proposition} \ 
Let $T\in(0,+\infty)$, and functions $u^1$ and $u^2\in C^1([0,T],L^\infty
$ 
$
(\mathbb R))$. 
Suppose that for any $t\in[0,T]$, 
the inequality 
\[u^1_t-\left(\mu*u^1-u^1+f(u^1)\right) 
\leq u^2_t-\left(\mu*u^2-u^2+f(u^2)\right)\] 
holds almost everywhere in $x$. Then, 
the inequality $u^1(T,x)\leq u^2(T,x)$ holds almost everywhere in $x$ 
if the inequality $u^1(0,x)\leq u^2(0,x)$ holds almost everywhere in $x$. 
\end{proposition}
{\bf Proof.} \ 
Put $K\in\mathbb R$ by 
\begin{equation}K:=1-\inf_{h>0,u\in\mathbb R}\frac{f(u+h)-f(u)}{h},\end{equation} 
and $v\in C^1([0,T],L^\infty(\mathbb R))$ by 
\begin{equation}v(t):=e^{Kt}(u^2-u^1)(t).\end{equation} 
Then, we have the ordinary differential equation 
\begin{equation}\frac{dv}{dt}=F(t,v)\end{equation} 
in $L^\infty(\mathbb R)$ with $v(0)=(u^2-u^1)(0)$ 
as we define a map 
$F:[0,T]\times L^\infty(\mathbb R)\rightarrow L^\infty(\mathbb R)$ by 
\[F(t,w):=\mu*w+(K-1)w
+e^{Kt}\left(f(u^1(t)+e^{-Kt}w)-f(u^1(t))\right)
+e^{Kt}a(t),\] 
where 
\[a:=\left(\frac{du^2}{dt}-\left(\mu*u^2-u^2+f(u^2)\right)\right)
-\left(\frac{du^1}{dt}-\left(\mu*u^1-u^1+f(u^1)\right)\right).\] 
For any $t\in[0,T]$, we see the inequality 
\begin{equation}a(t,x)\geq 0\end{equation}
almost everywhere in $x$.
Take the solution $\tilde{v}\in C^1([0,T],L^\infty(\mathbb R))$ to 
\begin{equation}\tilde{v}(t)=v(0)+\int_0^t\max\{F(s,\tilde{v}(s)),0\}ds.\end{equation} 
Then, for any $t\in[0,T]$, we have 
\begin{equation}\tilde{v}(t,x)\geq v(0,x)=(u^2-u^1)(0,x)\geq 0\end{equation}
almost everywhere in $x$. By using (5.1), (5.4) and (5.6), for any $t\in[0,T]$, 
we also have the inequality $F(t,\tilde{v}(t))\geq 0$ almost everywhere in $x$. 
Hence, from (5.5), $\tilde{v}(t)$ is the solution 
to the same ordinary differential equation (5.3) in $L^\infty(\mathbb R)$ 
as $v(t)$ with $\tilde{v}(0)=v(0)$. So, in virtue of (5.2) and (5.6), 
\[(u^2-u^1)(T,x)=e^{-KT}v(T,x)=e^{-KT}\tilde{v}(T,x)\geq 0\]
holds almost everywhere in $x$. 
\hfill 
$\blacksquare$ 

\vspace*{0.8em} 

\noindent
{\bf Proof of Proposition 13.}

The constants $0$ and $1$ are solutions to (4.1). So, by using Proposition 18, 
for any $u_0\in L^\infty(\mathbb R)$ with $0\leq u_0\leq 1$, 
there exists a solution $\{u(t)\}_{t\in[0,+\infty)}$ 
to (4.1) with $u(0)=u_0$ and $0\leq u(t)\leq 1$.

Let $u_0\in\mathcal M$. We take the solution $\{u(t)\}_{t\in[0,+\infty)}$ 
to (4.1) with $u(0)=u_0$. Let $t\in[0,+\infty)$ and $h\in[0,+\infty)$. 
Then, by Proposition 18, we see $u(t,x)\leq u(t,x+h)$ almost everywhere in $x$. 
We take a cut of function $\rho\in C^\infty(\mathbb R)$ with 
\[|x|\geq 1/2 \, \Longrightarrow \, \rho(x)=0,\]
\[|x|<1/2 \, \Longrightarrow \, \rho(x)>0\] 
and
\[\int_{x\in\mathbb R}\rho(x)dx=1.\] 
As we put  
\[v_n(x):=\int_{y\in \mathbb R}2^n\rho(2^n(x-y))u(t,y)dy\]
for $n\in\mathbb N$, we see $v_n(x)\leq v_n(x+h)$ for all $x\in\mathbb R$. 
Therefore, $v_n$ is a smooth and monotone nondecreasing function. 
By Helly's theorem, there exist a subsequence $n_k$ and $\psi\in\mathcal M$ 
such that $\lim_{k\rightarrow\infty}v_{n_k}(x)=\psi(x)$ almost everywhere in $x$. 
Then, $\|u(t,x)-\psi(x)\|_{L^1([-C,+C])}
\leq \lim_{k\rightarrow\infty}
(\|u(t,x)-v_{n_k}(x)\|_{L^1([-C,+C])}
+\|v_{n_k}(x)-\psi(x)\|_{L^1([-C,+C])})=0$ holds for all $C\in (0,+\infty)$. 
Hence, we obtain $\|u(t,x)-\psi(x)\|_{L^\infty(\mathbb R)}=0$. 
\hfill 
$\blacksquare$ 

\vspace*{0.4em}

\begin{proposition} \ 
Let a Borel-measure $\mu$ have $\lambda\in(0,+\infty)$ satisfying (4.2), 
and $T\in(0,+\infty)$. 
Suppose a sequence $\{u_n\}_{n=0}^\infty\subset C^1([0,T],L^\infty(\mathbb R))$ 
of solutions to (4.1) with $\sup_{n\in\mathbb N,x\in\mathbb R}|u_n(0,x)-u_0(0,x)|\leq 1$ 
satisfies \[\lim_{n\rightarrow\infty}\sup_{x\in\mathbb [-I,+I]}|u_n(0,x)-u_0(0,x)|=0\] 
for all $I\in(0,+\infty)$. Then, 
\[\lim_{n\rightarrow\infty}\sup_{t\in[0,T]}\|u_n(t,x)-u_0(t,x)\|_{L^\infty([-J,+J])}=0\] 
holds for all $J\in(0,+\infty)$. 
\end{proposition}
{\bf Proof.} \ 
Let $J\in(0,+\infty)$ and $\varepsilon\in(0,+\infty)$. 
We take $K\in[0,+\infty)$ such that  
\[K\geq\int_{y\in\mathbb R}e^{\lambda|y|}d\mu(y)
-1+\sup_{h>0,u\in\mathbb R}\frac{f(u+h)-f(u)}{h}.\] 
Put positive constants $\delta:=\min\{\varepsilon e^{-(KT+\lambda J)},1\}$ 
and $I:=\frac{1}{\lambda}\log(\frac{2}{\delta})$. 
Let $n\in\mathbb N$ be sufficiently large. Then, we have 
\begin{equation}\sup_{x\in\mathbb [-I,+I]}|u_n(0,x)-u_0(0,x)|\leq\delta.\end{equation} 
We consider the following two functions 
\[\underline{v}(t,x):=u_0(t,x)-e^{Kt}w(x)\]
and
\[\overline{v}(t,x):=u_0(t,x)+e^{Kt}w(x),\]
where $w(x):=\min\{\delta\frac{e^{\lambda x}+e^{-\lambda x}}{2},1\}$. 
We see 
\[(\mu*w)(x)\]
\[\leq \min\left\{\delta
\frac{
\left(\int_{y\in\mathbb R}e^{-\lambda y}d\mu(y)\right)e^{\lambda x}
+\left(\int_{y\in\mathbb R}e^{\lambda y}d\mu(y)\right)e^{-\lambda x}
}
{2},\mu(\mathbb R)\right\}\]
\[\leq \left(\int_{y\in\mathbb R}e^{\lambda |y|}d\mu(y)\right)w(x).\]
So, $\overline{v}$ is a super-solution to (4.1), because of 
\[\frac{d\overline v}{dt}-\left(\mu*\overline v-\overline v+f(\overline v)\right)\]
\[=(K+1)e^{Kt}w-\left(e^{Kt}(\mu*w)+\left(f(u_0+e^{Kt}w)-f(u_0)\right)\right)\geq 0\] 
almost everywhere in $x$. 
We can also see that $\underline{v}$ is a sub-solution. Because of $w(0)=\delta$, 
$w(\pm I)=1$ and (5.7), we get 
$\underline{v}(0,x)\leq u_n(0,x)\leq \overline{v}(0,x)$. Hence, by Proposition 18, 
$\underline{v}(t,x)\leq u_n(t,x)\leq \overline{v}(t,x)$ holds almost everywhere in $x$. 
So, we have $\|u_n(t,x)-u_0(t,x)\|_{L^\infty([-J,+J])}\leq e^{KT}w(\pm J)\leq \varepsilon$. 
\hfill 
$\blacksquare$

\vspace*{0.8em}

In virtue of Propositions 13, 18 and 19, if $\mu$ has a constant $\lambda\in(0,+\infty)$ 
satisfying (4.2), then  $Q^t \ (t\in(0,+\infty))$ satisfies Hypotheses 1 and $Q$ 5 
for the semiflow $Q=\{Q^t\}_{t\in[0,+\infty)}$ on the set $\mathcal M$ 
generated by (4.1). So, Theorems 6 and 7 can work for this semiflow. 

\vspace*{0.8em} 

\noindent
{\bf Proof of Theorem 14.} 

Put monotone nondecreasing functions 
$\varphi(x):=\max\{\alpha \in\mathbb R \, | \, 
\alpha \leq u(0,y)
$ 
$
\text{holds almost everywhere in } y\in (x,+\infty)\}$ 
and $\phi(x):=\lim_{h\downarrow+0}\varphi(x-h)$. Then, $\phi\in\mathcal M$, 
$\phi(-\infty)<1$ and $\phi(+\infty)=1$ hold.
We take a cut of function $\rho\in C^\infty(\mathbb R)$ with 
\[|x+1/2|\geq 1/2 \, \Longrightarrow \, \rho(x)=0,\]
\[|x+1/2|<1/2 \, \Longrightarrow \, \rho(x)>0\] 
and
\[\int_{x\in\mathbb R}\rho(x)dx=1.\] 
As we put  
\[v_n(x):=\int_{y\in \mathbb R}2^n\rho(2^n(x-y))u(0,y)dy\]
for $n\in\mathbb N$, we see $\phi\leq v_n$. Let $N\in\mathbb N$. 
Because of $\lim_{n\rightarrow\infty}\|v_{n}(x)-u(0,x)\|_{L^1([-N,+N])}=0$, 
there exists a subsequence $n_k$ such that $\lim_{k\rightarrow\infty}v_{n_k}(x)
$ 
$
=u(0,x)$ 
almost everywhere in $x\in[-N,+N]$. Therefore, we have 
$\phi(x)\leq u(0,x)$ almost everywhere in $x\in\mathbb R$. 
So, by Proposition 18, we obtain 
$Q^\tau[\phi](x-c\tau)\leq u(\tau,x-c\tau)=u(0,x)$ almost everywhere in $x$. 
Hence, because $Q^\tau[\phi](x-c\tau)\leq \varphi(x)$ holds, 
we get $Q^\tau[\phi](x-c\tau)\leq \phi(x)$. Therefore, by Theorem 6, 
there exists $\psi\in\mathcal M$ with $\psi(-\infty)=0$ and $\psi(+\infty)=1$
such that $Q^t[\psi](x-ct)\equiv\psi(x)$ holds for all $t\in[0,+\infty)$. 
\hfill 
$\blacksquare$ 

\vspace*{0.8em} 

\noindent
{\bf Proof of Theorem 15.} 

By Theorem 7,  there exists $c_*\in(-\infty,+\infty]$ such that the following holds: 
{\it Let $c\in \mathbb R$. Then, there exists $\psi\in\mathcal M$ 
with $\psi(-\infty)=0$ and $\psi(+\infty)=1$ 
such that $\{\psi(x+ct)\}_{t\in\mathbb R}$ is a solution to (4.1) 
if and only if $c\geq c_*$.}  

We show $c_*\not=+\infty$. 
Take $K\in[0,+\infty)$ such that  
\[K\geq\max\left\{\int_{y\in\mathbb R}e^{-\lambda y}d\mu(y),\mu(\mathbb R)\right\}
-1+\sup_{h>0}\frac{f(h)}{h}.\] 
As we put $\phi(x):=\min\{e^{\lambda x},1\}\in\mathcal M$, 
we see 
\[(\mu*\phi)(x)\leq\min
\left\{\left(\int_{y\in\mathbb R}e^{-\lambda y}d\mu(y)\right)e^{\lambda x},
\mu(\mathbb R)\right\}\]
\[\leq \max\left\{\int_{y\in\mathbb R}e^{-\lambda y}d\mu(y),\mu(\mathbb R)\right\}\phi(x).\]
So, $e^{Kt}\phi(x)$ is a super-solution to (4.1), because of 
\[e^{Kt}(\mu*\phi)-e^{Kt}w+f(e^{Kt}\phi)\leq Ke^{Kt}\phi.\] 
Hence, by Proposition 18, we obtain 
$Q^1[\phi](x)\leq e^K\phi(x)\leq e^{\lambda(x+\frac{K}{\lambda})}$, 
and $Q^1[\phi](x-\frac{K}{\lambda})\leq \phi(x)$. Therefore, from Theorem 6, 
there exists $\psi\in\mathcal M$ with $\psi(-\infty)=0$ and $\psi(+\infty)=1$
such that $Q^t[\psi](x-\frac{K}{\lambda}t)\equiv\psi(x)$ holds for all $t\in[0,+\infty)$. 
So, $c_*\leq \frac{K}{\lambda}$ holds. 
\hfill 
$\blacksquare$

\vspace*{0.8em} 

\noindent
{\bf Proof of Proposition 16.} 

Suppose $\psi$ is a continuous function. 
We take a interval $(a,b)\subset(0,1)$ such that 
\begin{equation}\inf_{u\in(a,b)}f_u(u)>1.\end{equation} 
Let $x\in\psi^{-1}((a,\frac{a+b}{2}))$ and $y\in\psi^{-1}((\frac{a+b}{2},b))$. 
Then, because of $x<y$ and (4.8), we have 
\[(\mu*\psi)(x)-\psi(x)+f(\psi(x))\]
\[\leq(\mu*\psi)(y)-\psi(x)+f(\psi(x))\]
\[<(\mu*\psi)(y)-\psi(y)+f(\psi(y)).\] 
It is a contradiction, as $\psi^{-1}((a,\frac{a+b}{2}))$ 
and $\psi^{-1}((\frac{a+b}{2},b))$ are open intervals. 
\hfill 
$\blacksquare$

\newpage

\noindent
{\bf Proof of Proposition 17.} 

We have $g(0)=1$ and $g^\prime(0)=-\int_{y\in\mathbb R}yd\mu(y)<0$ 
for $g(\zeta):=\int_{y\in\mathbb R}e^{-\zeta y}d\mu(y)
$ 
$
(\zeta\in[-\lambda,+\lambda])$. 
Hence, there exists $\xi\in(0,+\infty)$ with 
$\int_{y\in\mathbb R}e^{-\xi y}d\mu(y)<1$. 
Then, we take $\gamma\in(0,1-\int_{y\in\mathbb R}e^{-\xi y}d\mu(y))$. 

We consider the equation 
\begin{equation}u_t=\mu*u-u+\tilde f(u)\end{equation} 
instead of (4.1), 
where $\tilde f(u):=\gamma u$ in $(-\infty,0)$, $f(u)$ in $[0,1]$, 
$-2(u-1)$ in $(1,2)$ and $-u$ in $[2,+\infty)$. 
Also, we put $K:=\int_{y\in\mathbb R}e^{-\xi y}d\mu(y)-1+\gamma\in(-1,0)$. 
Then, we show that the function $\overline v(t,x):=2e^{K(t-1)}\min\{e^{\xi x},1\}$ 
is a super-solution to (5.9) on $t\in[0,1]$. 
For $x\in(-\infty,0]$, we can see 
$\mu*\overline v-\overline v+\tilde f(\overline v)
\leq \frac{d\overline v}{dt}$ 
from 
$(\mu*\overline v)(t,x)\leq 
2e^{K(t-1)}(\int_{y\in\mathbb R}e^{-\xi y}d\mu(y))e^{\xi x}$, 
$\tilde f(u)\leq \gamma u$ 
and $\frac{d\overline v}{dt}(t,x)=2Ke^{K(t-1)}e^{\xi x}$. 
For $x\in[0,+\infty)$, we can also see it 
from $(\mu*\overline v)(t,x)\leq 2e^{K(t-1)}$, 
$\tilde f(\overline v(t,x))=-\overline v(t,x)$ 
and $\frac{d\overline v}{dt}(t,x)=2Ke^{K(t-1)}$. 
Hence, by Proposition 18, we obtain 
$Q^1[\phi](x)\leq \overline v(1,x)\leq 2e^{\xi x}$, 
as we put $\phi(x):=\min\{2e^{\xi x-K},1\}\in\mathcal M$. 
So, because $Q^1[\phi](x-\frac{K}{\xi})\leq \phi(x)$ holds, 
by Theorem 6, there exists $\psi\in\mathcal M$ 
with $\psi(-\infty)=0$ and $\psi(+\infty)=1$ 
such that $\{\psi(x+\frac{K}{\xi}t)\}_{t\in\mathbb R}$ 
is a solution to (5.9). 
Because it is also one to (4.1), by Theorem 15, we have 
$c_*\leq \frac{K}{\xi}<0$. 

\hfill 
$\blacksquare$

\vspace*{3.2em} 

\noindent
{\bf A discrete Schrodinger model:}

We propose a discrete Schrodinger model that describes the measurement process. 
Let $X$ be the Hilbert space of state vectors 
of a one-particle system of spin $\displaystyle\frac{1}{2}$ 
on the one-dimensional discrete grid $\mathbb Z$. 
That is, $X=L^2(\, {\mathbb Z}\, ; \, {\mathbb C}^2\, )$.
Let natural numbers $L_0$ and $N_0$ be large. 
Put $I=\{ \, n\in {\mathbb Z} \, |\,  L_0\leq|n|\leq L_0+N_0\, \}$. 
We assume that $I$ is {\bf the place where two detectors exist}. 
Suppose that \underline{for $n\in I$}, $V(n)$ is a $2\times 2$ Hermitian {\bf random matrix}. 
Suppose that for $n\in I$, if $U$ is a $2\times 2$ unitary matrix, 
then the distributions of  $(U^{-1})(V(n))U$ and $V(n)$ are the same.
Suppose that for $n\in I$, the density function of $V(n)$,  
roughly speaking, is smooth and compact supported. 
Suppose that $\{V(n)\}_{n\in I}$ is independent and identically distributed. 
Suppose that \underline{for $n\not\in I$}, 
$V(n)$ is the $2\times 2$ {\bf zero matrix}. 
Then, we propose a discrete Schrodinger model 
on the Hilbert space $X\ (\, =L^2(\, {\mathbb Z}\, ; \, {\mathbb C}^2\, )\, )$
$$\sqrt{-1}\, \frac{du}{dt}\, (n)\, 
=\, -\frac{1}{2m}\, (\, (u(n-1)-u(n))\, +\, (u(n+1)-u(n))\, )\, 
+\, V(n)\, u(n)\, .$$

\newpage

\noindent Acknowledgments. \ 
I thank Professor Hiroshi Matano, Mr. Xiaotao Lin 
and Mr. Masahiko Shimojo for their discussion. It was partially supported 
by Grant-in-Aid for Scientific Research (No.19740092) 
from Ministry of Education, Culture, Sports, Science 
and Technology, Japan. 


\[ \begin{array}{c} \mbox{R\scriptsize EFERENCES}  \end{array} \]

[1] P. W. Bates, P. C. Fife, X. Ren and X. Wang, Traveling waves in a convolution model 
for phase transitions, {\it Arch. Rational Mech. Anal.}, 138 (1997), 105-136. 

[2] J. Carr and A. Chmaj, Uniqueness of travelling waves for nonlocal monostable equations, 
{\it Proc. Amer. Math. Soc.}, 132 (2004), 2433-2439. 

[3] X. Chen, Existence, uniqueness, and asymptotic stability of traveling waves 
in nonlocal evolution equations, {\it Adv. Differential Equations}, 2 (1997), 125-160. 

[4] X. Chen and J.-S. Guo, Uniqueness and existence of traveling waves for discrete quasilinear monostable dynamics, {\it Math. Ann.}, 326 (2003), 123-146. 

[5] J. Coville and L. Dupaigne, On a non-local equation arising in population dynamics, 
{\it Proc. Roy. Soc. Edinburgh A}, 137 (2007), 727-755.

[6] R. A. Fisher, {\it The Genetical Theory of Natural Selection}, Clarendon Press, Oxford, 1930. 

[7] R. A. Fisher, The wave of advance of advantageous genes, {\it Ann. Eugenics}, 7 (1937), 335-369. 

[8] B. H. Gilding and R. Kersner, {\it Travelling Waves in Nonlinear 
Diffusion-Convection Reaction}, Birkh\"auser, Basel, 2004. 

[9] J.-S. Guo and F. Hamel, Front propagation for discrete periodic monostable equations, 
{\it Math. Ann.}, 335 (2006), 489-525. 

[10] J.-S. Guo and Y. Morita, Entire solutions of reaction-diffusion equations 
and an application to discrete diffusive equations, {\it Discrete Contin. Dynam. Systems}, 
12 (2005), 193-212. 

[11] F. Hamel and N. Nadirashvili, Entire solutions of the KPP equation, 
{\it Comm. Pure Appl. Math.}, 52 (1999), 1255-1276. 

[12] Y. Hosono, The minimal speed for a diffusive Lotka-Volterra model, {\it Bull. Math. Biol.}, 
60 (1998), 435-448. 

[13] Y. Kan-on, Fisher wave fronts for the Lotka-Volterra competition model with diffusion, 
{\it Nonlinear Anal.}, 28 (1997), 145-164. 

[14] A. N. Kolmogorov, I. G. Petrovsky and N. S. Piskunov, \'Etude de l'\'equation de la difusion 
avec croissance de la quantit\'e de mati\`ere et son application \`a un probl\`eme biologique, 
{\it Bull. Univ. Moskov. Ser. Internat. A}, 1 (1937), 1-25. 

[15] B. Li, H. F. Weinberger and M. A. Lewis, Spreading speeds as slowest wave speeds 
for cooperative systems, {\it Math. Biosci.}, 196 (2005), 82-98. 

[16] X. Liang and X.-Q. Zhao, Asymptotic speeds of spread and traveling waves 
for monotone semiflows with applications, {\it Comm. Pure Appl. Math.}, 60 (2007), 1-40. 

[17] H. Okamoto and M. Shoji, {\it The Mathematical Theory of Permanent Progressive Water-Waves}, 
World Scientific Publishing Co., River Edge, 2001. 

[18] K. Schumacher, Travelling-front solutions for integro-differential equations. I, 
{\it J. Reine Angew. Math.}, 316 (1980), 54-70. 

[19] K. Schumacher, Travelling-front solutions for integrodifferential equations II, 
{\it Biological Growth and Spread}, pp. 296-309, Springer, Berlin-New York, 1980. 

[20] K. Uchiyama, The behavior of solutions of some nonlinear diffusion equations for large time, 
{\it J. Math. Kyoto Univ.}, 18 (1978), 453-508. 

[21] A. I. Volpert, V. A. Volpert and V. A. Volpert, {\it Traveling wave solutions 
of parabolic systems}, American Mathematical Society, Providence, 1994. 

[22] H. F. Weinberger, Long-time behavior of a class of biological models, 
{\it SIAM J. Math. Anal.}, 13 (1982), 353-396. 

[23] H. F. Weinberger, On spreading speeds and traveling waves for growth and migration models 
in a periodic habitat, {\it J. Math. Biol.}, 45 (2002), 511-548. 

[24] H. F. Weinberger, M. A. Lewis and B. Li, Analysis of linear determinacy for spread 
in cooperative models, {\it J. Math. Biol.}, 45 (2002), 183-218.

\end{document}